\documentclass[11pt]{article}
\usepackage{amsmath}
\usepackage{amssymb}
\usepackage{eucal}
\usepackage{stmaryrd}

\begin{document}

\baselineskip 16pt

\title{On one generalization of finite $\frak{U}$-critical groups}

\author{Vladimir N. Semenchuk\\
{\small 
Department of Mathematics,  Francisk Skorina Gomel State University,}\\
{\small Gomel 246019, Belarus}\\
{\small E-mail:  kolenchukova@gsu.by}\\ \\
{ Alexander  N. Skiba 
}\\
{\small Department of Mathematics,  Francisk Skorina Gomel State University,}\\
{\small Gomel 246019, Belarus}\\
{\small E-mail: alexander.skiba49@gmail.com}}

\date{}
\maketitle

\begin{abstract} A  proper  subgroup $H$ of  a group $G$ is said to be: 
\emph{$\Bbb{P}$-subnormal} in $G$  if  there exists a chain of subgroups 
$H=H_0 < H_1< \cdots < H_{n}=G$ such that 
  $|H_{i}:H_{i-1}|$ is a prime for $i=1,\ldots
,n$; \emph{$\Bbb{P}$-abnormal} in $G$ if for every two subgroups $K\leq L$
 of $G$, where $H\leq K$,  $|L:K|$ is not a prime. 
 In this paper we describe finite groups in which every non-identity subgroup is either 
$\Bbb{P}$-subnormal  or $\Bbb{P}$-abnormal.

\end{abstract}

\footnotetext{Keywords: finite group, $\Bbb{P}$-subnormal subgroup,
 $\Bbb{P}$-abnormal subgroup,  $\frak{U}$-critical group,   Gasch\"utz subgroup, Hall subgroup
}

\footnotetext{Mathematics Subject Classification (2010): 20D10,
20D15, 20D20}
\let\thefootnote\thefootnoteorig

\section{Introduction}

Throughout this paper, all groups are finite, $G$ denotes a finite
group and $p$ is a prime. We use $\frak{N}$   and $\frak{U}$ to denote the classes 
 of all nilpotent and of all supersoluble groups, respectively.
A sugroup  $H$ of $G$ is said to be a  \emph{Gasch\"utz
} subgroup of $G$ (Shemetkov \cite[p. 170]{Shem}) if 
 $H$ is  supersoluble and   $|L:K|$ is not a prime whenever 
$H\leq K \leq L \leq G$.

Let $\frak{F}$ be a class of groups. If $1\in {\frak{F}}$, then we
 write $G^{\frak{F}}$ to denote the
intersection of all normal subgroups $N$ of $G$ with  $G/N\in {\frak{F}}$. The
 class $\frak{F}$ is said to be a \emph{formation} if either
${\frak{F}}= \varnothing $ or $1\in {\frak{F}}$ and every
homomorphic image of $G/G^{\frak{F}}$  belongs to $ {\frak{F}}$ for any group $G$.
 The formation
${\frak{F}}$ is said to be: \emph{saturated}  if $G\in {\frak{F}}$
whenever $G^{\frak{F}}\leq \Phi (G)$; \emph{
 hereditary} if $H\in {\frak{F}}$ whenever $ G \in {\frak{F}}$
   and $H$ is a  subgroup of $G$.

  A group $G$ is said to be   \emph{${\frak{F}}$-critical } if $G$ is not
 in ${\frak{F}}$ but all proper subgroups of $G$ 
are in ${\frak{F}}$ \cite[p. 517]{DH}.
 An ${\mathfrak N}$-critical group is also   called   a \emph{Schmidt group}.

A proper subgroup $H$ of  $G$ is said to be: 
\emph{$\frak{F}$-subnormal} in $G$ if  there exists a chain of subgroups 
$H=H_0 <  H_1 <  \cdots <  H_{n}=G$ such that  $H_{i-1}$ is a maximal 
subgroup of $H_{i}$
 and $H_{i}/(H_{i-1})_{H_i}\in \frak{F}$ for all $i=1,\ldots
,n$; \emph{$\frak{F}$-abnormal} in $G$ if   $L/K_{L}\not \in \frak{F}$ 
whenever    $H\leq K  < L \leq G$ and  $K$ is a maximal subgroup of $L$.
 A group $G\not \in \frak{F}$ is said to be an \emph{$E_{\frak{F}}$-group} \cite{2} if 
every non-identity  
subgroup of $G$  is  either $\frak{F}$-subnormal   
or $\frak{F}$-abnormal  in $G$.

In \cite{1}, Fattahi described  groups in which every 
subgroup is either normal or abnormal.   As a generalization of this result, 
Ebert and  Bauman  classified the $E_{\frak{F}}$-groups in the case when 
${\frak{F}}={\frak{N}}$ (in this case $G$ is a group in which
 every subgroup is either subnormal or abnormal), and in the case when ${\frak{F}}$
 is the class of all soluble $p$-nilpotent groups, for  odd prime $p$
  \cite{2}. In the future, the $E_{\frak{F}}$-groups were  studied for some other
${\frak{F}}$ (see for example \cite{3, I, II, III}). 
Nevertheless, it should be noted that
 a complete description of the $E_{\frak{F}}$-groups  was obtained 
only for  such cases ${\frak{F}}$  when every ${\mathfrak F}$-critical group 
  is  a  Schmidt group 
 \cite{1, I, II, III}). Thus, for example 
 in the  case, where  $\frak{F}= {\frak{U}}$, the structure  of
  $E_{\frak{F}}$-groups has not been known since the methods in \cite{1, 
3, I, II, III} could   not be used  in the analysis of this case.

Note, in passing,  that if $G$ is soluble and $H$ is a subgroup of $G$, then $H$: is 
$\frak{U}$-subnormal in $G$ if and only   if  there exists a chain of subgroups 
$H=H_0 < H_1< \cdots < H_{n}=G$ such that 
  $|H_{i}:H_{i-1}|$ is a prime for $i=1,\ldots
,n$; $\frak{U}$-abnormal in $G$ if   and only if  $|L:K|$ is not a prime whenever 
$H\leq K \leq L \leq G$.

If $G$ is supersoluble, then clearly every subgroup of $G$ is ${\mathfrak U}$-subnormal 
in $G$.  
A full  description  of  $E_{\frak{U}}$-groups,  for the non-supersoluble case,
 gives the following our result.

{\bf  Theorem A. }   {\sl Let $G$ be an  $E_{\frak{U}}$-group 
 and     $D=G^{\frak{U}}$
 the supersoluble   residual of $G$.
Then  $G=D\rtimes H$, where: } 

(i) {\sl     $H$  is a Hall Gasch\"utz subgroup   of $G$. Hence 
 if $H$ is nilpotent,
 then it is  a Carter subgroup of $G$.  }

(ii) {\sl Every chief factor of $G$ below $D$ is non-cyclic. Hence $H$ is a
 supersoluble normalizer ($\frak{U}$-normalizer, in other words) of $G$
 in the sence of} \cite{Car}.

(iii) {\sl $|G:DG'|$ is a prime power number.} 

(iv) {\sl If $H$ is not a  cyclic group of prime power order $p^{n}$, where $n > 1$,
 then $D$ is nilpotent. }

(v) {\sl $H\Phi (G)/\Phi (G)$ is
 either a Miller-Moreno group     or an abelian
 group of prime power order.}

(vi) {\sl Every proper subgroup of $G$ containing $D$ is supersoluble.

Conversely, any group satisfying the above conditions is an $E_{\frak{U}}$-group. }

{\bf Corollary 1.1.}  {\sl Let $G=D\rtimes H$ be an   $E_{\frak{U}}$-group,
 where  
      $D=G^{\frak{U}}$.  If $H$ is nilpotent, then 
 it is a system   normalizer of $G$.  }

From the  describtion  of ${\frak{U}}$-critical groups $G$  \cite{hu, D}  it follows 
that every  subgroup of $G$  contaning $\Phi (G) \cap G^{\frak{U}}$ is
 either $\frak{U}$-subnormal    or $\frak{U}$-abnormal  in $G$ (see Lemma 2.6 below).  
Another application of Theorem A is
 the following result, which classifies  all  groups with such a property.

{\bf  Theorem B. }   {\sl Let $G$ be a non-supersoluble group  
 and $\Phi =\Phi (G) \cap G^{\frak{U}}$.
  Then every non-identity  subgroup of   $G$ contaning $\Phi$ is
 either $\frak{U}$-subnormal   
or $\frak{U}$-abnormal  in $G$     if
 and only if  $G=D\rtimes H$ is a soluble group, where  $H$ is a  
  Hall subgroup of $G$ such that  $H\Phi /\Phi$ is a Gasch\"utz subgroup
 of $G$, and  with respect to $G$  Assertions (iii)--(vi) in Theorem A hold.  }

All unexplained notation and terminology are standard. The reader is
referred to \cite{Shem-Sk}, \cite{DH},  \cite{Guo},  or  \cite{BE}  if necessary.

\section{Preliminaries}

The following lemma collects some well-known properties of $\frak{F}$-subnormal 
subgroups which will be used in our proofs.

{\bf Lemma 2.1.} {\sl Let $\mathfrak F$   be a hereditary saturated formation,  
 $H$ and $K$  the subgroups of $G$ and     $H$
is $\mathfrak F$-subnormal in $G$. }

(1) {\sl $H\cap K$
is  $\mathfrak F$-subnormal in
$K$} \cite[6.1.7(2)]{BE}.

(2) {\sl If $N$ is a normal subgroup in $G$, then $HN/N$ is
$\mathfrak F$-subnormal in $G/N$.    } \cite[6.1.6(3)]{BE}.

(3) {\sl If  $K$ is an $\mathfrak F$-subnormal subgroup of  $H$,
then $K$ is $\mathfrak F$-subnormal in $G$} \cite[6.1.6(1)]{BE}.

(4) {\sl If $G^{\mathfrak F}\leq K$, then $K$
is $\mathfrak F$-subnormal in $G$} \cite[6.1.7(1)]{BE}.

(5) {\sl If $K\leq H$ and $H\in \mathfrak F$, then $K$
is $\mathfrak F$-subnormal in $G$.}

A minimal normal subgroup $R$   of $G$ is 
called \emph{$\mathfrak{F}$-central in $G$} provided $R\rtimes
(G/C_{G}(R))\in \mathfrak{F}$, otherwise it is called  
\emph{$\mathfrak{F}$-eccentric} in $G$.

{\bf Lemma 2.2.} {\sl Let $ \mathfrak F$  be a formation and  $M$  a maximal subgroup of $G$.
Let  $R$  be 
 a minimal normal subgroup of $G$ such that $MR=G$. Then $G/M_{G}\in \mathfrak F$ if and only if $R$ is 
 $ \mathfrak F$-central in $G$.  }

{\bf Proof. }   In view of the $G$-isomorphism $R\simeq RM_{G}/M_{G}$ we 
can assume without loss of generality that $M_{G}=1$. Let 
$C=C_{G}(R)$. 

 If $G\simeq G/M_{G} \in  
\mathfrak F$, then $R$ is   $\mathfrak F$-central in $G$ by the  
 Barnes-Kegel's  Theorem  \cite[IV, 1.5]{DH}. 
 Now let    $R\rtimes (G/C)\in \mathfrak{F}$. First assume that $R$ 
is non-abelian. If $R$ is the unique minimal normal subgroup of $G$, then 
$C=1$ and so $G\in \mathfrak F$. Now let $G$ have a minimal normal subgroup 
$L\ne R$. Then, since  $M_{G}=1$, $G=R\rtimes M=L\rtimes M$ and $C=L$ by \cite[A, 15.2]{DH}. 
 Hence $M\simeq G/L\simeq G/R\in \mathfrak{F}$ and so $G\in \mathfrak{F}$ 
since $\mathfrak{F}$ is a formation. 

Finally, if $R$ is an abelian $p$-group, then $C=R$   by \cite[A, 15.2]{DH} 
and  so 
  $G\simeq G/M_{G}\simeq R\rtimes (G/R)\in  \mathfrak{F}$.    The lemma is 
proved.

 {\bf Lemma 2.3} (See Lemma 2.15 in \cite{5}). {\sl Let $E$ be a  normal non-identity 
quasinilpotent subgroup of $G$.
  If $\Phi (G)\cap E=1$, then $E$
 is the direct product of some minimal normal subgroups of $G$.}

{\bf Lemma 2.4.} {\sl Let  $\frak{F}$ be a non-empty hereditary saturated
 formation, $G$  an  $E_{\frak{F}}$-group and $D=G^{\frak{F}}$. }

 (i) {\sl Every  $\frak{F}$-subnormal subgroup of $G$ belongs to   
$\frak{F}$. }

(ii) {\sl  $F^{*}(G)\leq  D\Phi(G)$.}

{\bf Proof. }  (i)  Let  $H$ be any $\frak{F}$-subnormal subgroup of $G$ 
and $K$  a maximal subgroup of  $H$. Then  $K$ is not  
$\frak{F}$-abnormal, so it is $\frak{F}$-subnormal in $G$ by hypothesis. Hence $K$ is 
 $\frak{F}$-subnormal in  $H$ by Lemma 2.1.(1), that is, $H/K_{H}\in 
\frak{F}$.  
 Therefore, since  $\frak{F}$ is a saturated 
formation,  $H\in\frak{F}$.

(ii) Without loss of generality we can assume that    $\Phi(G)=1$.  In this 
case $F^{*}(G) =N_{1}\times \cdots \times N_{t}$ for
 some minimal normal subgroups  $N_{1}, \ldots , N_{t}$ of 
  $G$ by Lemma 2.3. Let  $N=N_{i}$ and  $M$  be a maximal subgroup of $G$ such that 
 $G=MN$. Assume that  $M$ is $\frak{F}$-subnormal in $G$. Then  $D\leq M_{G}$, so 
 $N$ is  $\frak{F}$-central in $G$ by Lemma 2.2. On 
the other hand,  Assertion (i) implies that 
 $M\in\frak{F}$.  Thus $G/N\simeq  M/M\cap N \in  \frak{F}$ and so
 $G\simeq G/N\cap M_{G}\in \frak{F}$. This 
contraiction shows that 
  $M$ is $\frak{F}$-abnormal in  $G$, so $G/M_{G}  \not \in \frak{F}$.
   Hence $N\leq  D$ by Lemma 2.2. 
Therefore $F^{*}(G)\leq D$. The lemma is proved.

{\bf Lemma 2.5.} {\sl Let  $\mathfrak F$  be a non-empty formation,  
$G$  an 
$\mathfrak{F}$-critical  group and  $D=G^{\frak F}$.
 }

(i) {\sl If $G$ is soluble , then 
$D$  is a $p$-group for some prime $p$.}

(ii) {\sl  If $\mathfrak F$  is saturated and $D$ is soluble, then the 
following statements hold:}

(a) {\sl $D$  is a $p$-group for some prime $p$.}

(b) {\sl $D/\Phi(D)$~ is a chief factor of $G$.}

{\bf Proof. } (i)   Since $G$ is soluble,  $\Phi(G) < F(G)$. Hence for some prime $p$ 
we have  $O_{p}(G) \nleq \Phi (G)$. Let $M$ be a maximal subgroup of $G$ such that 
 $G=O_{p}(G) M$. Then $G/O_{p}(G) =O_{p}(G) M/O_{p}(G) \simeq 
 M/M\cap O_{p}(G) \in \mathfrak F$ since $G$ is
 an $\mathfrak{F}$-critical group. Thus $D\leq O_{p}(G) $.  

(ii) See  Theorem 24.2 in \cite[V]{Shem} or  \cite[VII, 6.18]{DH}.  The lemma is proved.

{\bf Lemma 2.6.} {\sl Let  $\mathfrak F$  be a  hereditary
 saturated formation and    $G$
  an $\mathfrak{F}$-critical soluble group. 
 Then every
   subgroup of $G$ contaning $\Phi (G) \cap G^{\mathfrak F}$ is either  
$\frak{F} $-subnormal or $\frak{F} $-abnormal in $G$.}
                                                                    
{\bf Proof. }  It is enough to consider the case when $\Phi (G) \cap G^{\mathfrak F}=1$.
By Lemma 2.5,   $D$ is a minimal 
normal subgroup of $G$. Let $A$ be any non-identity subgroup of $G$. First 
assume that  $DA < G$. Then $DA\in \mathfrak{F}$, and $DA$
 is $\mathfrak{F}$-subnormal in $G$ by  Lemma  2.1(4). Hence $A$ 
 is $\mathfrak{F}$-subnormal in $G$ by Lemma  2.1(3). Now 
assume that $DA=G$. Then $A$ is a maximal subgroup of $G$, so  $A$  is 
$\mathfrak{F}$-abnormal in $G$.  The lemma is proved.

{\bf Lemma 2.7} (Friesen  \cite[4, 3.4]{We}).  {\sl If $G=AB$, where $A$ and $B$ are normal
 supersoluble subgroups of $G$ and  $(|G:A|, |G:B|)=1$, then $G$ is supersoluble.}

We shall need the following special case of Theorem C in \cite{5}.

 {\bf Lemma 2.8. }
  {\sl Let $\cal F $ be a  hereditary saturated
formation containing all nilpotent groups and  $E$ a normal subgroup
of $G$. If $E/E\cap \Phi (G) \in {\cal F}$,
 then $E \in {\cal F}$.}

A   subgroup $H$ of  $G$ is said to be: 
\emph{$\Bbb{P}$-subnormal} in $G$ \cite{11, 22} if  there exists a chain of subgroups 
$H=H_0 < H_1< \cdots < H_{n}=G$ such that 
  $|H_{i}:H_{i-1}|$ is a prime for $i=1,\ldots
,n$; \emph{$\Bbb{P}$-abnormal} in $G$ if $|L:K|$ is not a prime whenever 
$H\leq K \leq L \leq G$.  We say that $H$ satisfies 
   the \emph{$\Bbb{P}$-property} in $G$  if $H$ is     either $\Bbb{P}$-subnormal
 or $\Bbb{P}$-abnormal    in $G$.

 {\bf Lemma 2.9.}   (i)  {\sl If every  non-identity subgroup of $G$ of prime order satisfies  
the $\Bbb{P}$-property in $G$, then $G$ is not a simple non-abelian group.}

(ii) {\sl If every non-identity cyclic  subgroup of $G$ of prime power order 
satisfies  
the $\Bbb{P}$-property in $G$, then $G$ is soluble. }

 {\bf Proof.} (i) Suppose that this is false and let $p$ be
 the smallest prime dividing $|G|$. Then
 a Sylow $p$-subgroup $P$ of $G$ is not cyclic. Let $H$ be a subgroup of 
order $p$ in  $P$. Then $H  < P$, so by hypothesis, $G$ has a  maximal 
subgroup $M$ such that $H\leq  M$   and   $|G:M|=q$ for some prime $q$. 
                  Since $G$ is a simple non-abelian group, $M_G=1$ and 
by considering the permutation
representation of $G$ on the right cosets of $H$, we see that $G$ is isomorphic to some
subgroup of the symmetric group $S_q$ of degree $q$.
Hence   $q$ is the 
largest prime divisor of $|G|$   and  $|Q|=q$, where $Q$ is  a Sylow $q$-subgroup $Q$
 of $G$.  It follows that   $q\ne p$.  It is clear that $G$ is not 
$q$-nilpotent, so it has a $q$-closed Schmidt subgroup $H$ such that 
$Q\leq H$ by \cite[IV, 5.4]{hupp}. Since $Q$ is normal in $H$, it is 
$P$-subnormal  in $G$ by hypothesis. Hence $G$  has a  maximal subgroup $T$ 
such that $Q\leq T$ and $|G:T|=r$ is a prime. But then $r$ is the largest 
prime dividing $|G|$ and  so $r=q$, a 
contradiction. Hence we have (i). 

(ii) Since  the hypothesis clearly holds for every quotient of $G$ and every normal subgroup of
$G$, this  assertion is a corollary of Assertion (i). The lemma is proved.

\section{Proofs of Theorems A and B}

{\bf Proof of Theorem A.}   {\sl Necessity.} 
  Suppose that this is false and let $G$ be a counterexample 
 of   minimal order.  Let $\pi =\pi (D)$.

(1) {\sl The hypothesis holds  on  $G/R$ for every
  normal subgroup $R$ of $G$ not containing $D$.}

First note that $G/R\not \in  
\frak{U}$ since $D\nleq R$.  Therefore this claim is a corollary of Lemma 
2.1(2).

(2) {\sl Every subgroup $E$ of $G$ containing $D$ is supersoluble. Hence $G$ 
is soluble.}

First note that the hypothesis  holds for $D$, so $D$ is soluble by Lemma  
2.9. On the other hand,    $E$  is
  $\frak{U}$-subnormal in $G$ by Lemma 2.1(4) and  so $E$ is supersoluble by   
Lemma 2.4(i). Hence we have (2).

(3) {\sl $D$ is a Hall subgroup of $G$. }

Suppose that this is false and let $P$ be a Sylow $p$-subgroup of
$D$ such that $1 < P < G_{p}$, where $G_{p}$ is a Sylow $p$-subgroup of 
$G$. Then $|G_{p}| > p$.  Let  $R$ be a minimal normal subgroup of
$G$.

(a)  {\sl   $R$ is a $p$-group. Hence $O_{p'}(G)=1$}.

Since $G$ is soluble by Claim (2),    $R$ is a $q$-group  for some
prime $q$.  Moreover,
 $DR/R= (G/R)^{\frak{U}}$ is a Hall subgroup of $G$ by the choice 
of $G$ since the hypothesis holds for $G/R$ by Claim (1). Therefore every Sylow
 $r$-subgroup of $D$, where $r\ne q$,  is a Sylow subgroup of $G$. 
 Hence  $q=p$ and so  $O_{p'}(G)=1$. 

(b)  {\sl    If $R\leq D$, then $R$ is a  Sylow $p$-subgroup of $D$.
 If $R\cap D=1$, then $G_{p}=R\rtimes P$.}

Assume $R\leq D$. Then $R\leq P$ and $P/R$ is a Sylow $p$-subgroup of $D/R$.
 If $P/R\ne 1$, then Claim (1) and the choice of $G$ imply that 
$P/R=G_{p}/R$ and  so $P=G_{p}$. This contradiction shows that $P=R$
 is a Sylow $p$-subgroup of $D$.

Now  assume that $R\nleq  D$. Then $RP/R=G_{p}/R$  since   $DR/R=(G/R)^{\frak{U}}$ is a Hall 
subgroup of $G/R$ by Claim (1), so $G_{p}=R\rtimes P$ since $R\cap D=1$. 

(c) {\sl $R\nleq \Phi (G)$. Hence  $\Phi (G)=1$.
}

Assume that  $R\leq \Phi (G)$. If $R\cap D=1$, then $G_{p}=R\rtimes P$ by 
Claim (b) and so  $R\nleq \Phi (G)$ by the 
Gasch\"utz Theorem \cite[I,  17.4]{hupp}. This contradiction shows that 
$R\leq D$, so $R$ is the Sylow $p$-subgroup of $D$ by Claim (b). Hence $D/R$ is a $p'$-group,
 so a $p$-complement $S$ of $D$ is normal in $G$ by Lemma 2.8. But then 
$S\leq  O_{p'}(G)=1$. Hence $D=R$ and so $G$     is  supersoluble since 
$D=G^{\mathfrak{U}}$, a contrdiction.  Thus we have (c).

(d) {\sl $G_{p}$ is normal in $G$.}

Let $E$ be any normal maximal subgroup of $G$ containing $D$ with $|G:E|=q$.  Then 
$O_{p'}(E)\leq O_{p'}(G)=1$, so $p$ is the largest prime dividing $|E|$ since $E$
 is supersoluble by Claim (2). If $q\ne p$, then 
$G_{p}\leq E$ and so $G_{p}$ is normal in $G$ since in this case $G_{p}$ is a 
characteristic subgroup of $E$.

Finally, assume that  $q=p$. Then $p$ is the largest prime dividing $|G|$ and
so  $DG_{p}$ is normal in $G$ since $G/D=G/G^{\frak{U}}$ is supersoluble.
  If $DG_{p}\ne G$, we can get as above that 
$G_{p}$ is normal in $G$. Now assume that   $DG_{p}= G$.
 Since $R\nleq \Phi (G)$ by Claim (c), it has 
a complement in $G$ and so $R$ has a complement $V$ in $G_{p}$. It is 
clear that $V$ is not $\frak{U}$-abnormal in $G$, so for some maximal 
$\frak{U}$-subnormal subgroup $M$ of $G$ we have  $V\leq M$,  which implies 
that $G=DV\leq M$. This contradiction  shows that the case under 
consideration is impossible. Hence $G_{p}$ is normal in $G$

{\sl Final contradiction for (3).}  In view of Claim (d), $\Phi (G_{p})\leq \Phi (G)=1$.
 Therefore, by the Maschke's Theorem,  $G$
 has a minimal normal subgroup $L\nleq \Phi (G)$
 such that $L\leq G_{p}$  and  $L\nleq  D$. Then $|L|=p$ and for some maximal subgroup $M$ of 
$G$ we have $G=LM$. Hence $G/M_{G}$ is supersoluble by Lemma 2.3, which implies that
 $D\leq M$ and so $M$  
 is supersoluble by Claim (2). But then  $G$ is supersoluble, a contradiction. Hence we have (3).   

(4) {\sl $F(G)\leq D\Phi(G)$} (This directly follows from Lemma 2.4(ii)).

(5) {\sl $|G:DG'|$ is a prime power number} (Since $G$ is
 not supersoluble, this directly follows
from Claim (2) and Lemma 2.7).

(6) {\sl If $\Phi (G)=1$, then $O_{\pi'}(G)=1$. }

Assume that $O_{\pi'}(G)\ne 1$ and let $R$ be a minimal normal subgroup of $G$ contained in 
$O_{\pi'}(G)$. Then  $R\nleq  D$, so in view of the 
$G$-isomorphism $RD/D\simeq R$, $R$ is cyclic.
  Since $\Phi (G)=1$,
 for some maximal subgroup $M$ of $G$ we have   $RM=G$  and so $G/M_{G}$ is supersoluble
 by Lemma 2.3.  It follows that $D\leq M$ and hence $M$ is supersoluble by 
Claim (2).   But then   $G=RM$ is supersoluble. This 
contradiction shows that  we have (6).

(7) {\sl If $H$ is a complement to $D$ in $G$, then $H\Phi (G)/\Phi (G)$ is
 either a Miller-Moreno group or an abelian group of prime power order.}

Without loss of generality we can assume that    $\Phi(G)=1$. 
First we shall show that every proper subgroup $A$  of $H$ is abelian.  
Let $C=C_{G}(F(G))$. Then $C\leq F(G)$ since $G$ is soluble. On the other hand, Claim (4)
 implies that $F(G)\leq D$, so $C\leq D$. It follows that $F(DA)=F(G)$, so 
$AF(G)/F(G)\simeq A$ is abelian since $DA$ is supersoluble by Claim (2). 
Therefore, if $H$ is not abelian, then $H$ is a Miller-Moreno group.

Finally, suppose  that $H$ is abelian. Then $G'\leq D$ by Claim (3), so
   Claim (5) implies that
 $|G:DG'|=|G:D|=|H|$ is a prime power number.

(8) {\sl $H$  is a  Gasch\"utz subgroup   of $G$. }

Since $D=G^{\frak{U}}$ and $G=D\rtimes H$, $H$ is supersoluble. It is 
clear also that $H$ is not $\frak{U}$-subnormal in $G$. Hence $H$ 
is $\frak{U}$-abnormal in $G$ by hypothesis.  Therefore  $H$ 
 is a Gasch\"utz subgroup   of   
$G$ since $G$ is soluble by Claim (2).

(9) {\sl If $H$ is not a cyclic group of prime power order $q^{n}$, where $n > 1$,
  then $D$ is nilpotent. }
          
Suppose that this is false and let  $R\leq O_{p}(G)$ be a minimal normal subgroup of
$G$. 

(*) {\sl $|H|$ is not a prime}.

Indeed,  assume that $H=\langle a\rangle $, where $|a|$ is a prime. 
 Since   $H$ is a Gasch\"utz subgroup of $G$, $N_{G}(H)=H$ and 
hence $a$ induces a regular automorphism on $D$. Hence $D$ is nilpotent by 
the Thompson's theorem \cite[V, 8.14]{HupBl}, a contradiction.  Hence we 
have (*).

(**) {\sl If $R\leq D$ or $R\leq \Phi (G)$, where $R$ is a minimal normal subgroup of $G$,
 then  $DR/R$ is nilpotent. Hence $\Phi (G)=1$, $R=F(D)=C_{D}(R)$ is the unique
 minimal normal subgroup of $G$ contained in $D$ and $R$ is the Sylow $p$-subgroup of $G$
 for some prime $p$.
}

The choice of $G$ and Claims (1)  and (*) imply that in order to prove that $DR/R$ is nilpotent,
 it is enough to show that $HR/R$ is not a cyclic group of order $q^{n}$, where
 $n > 1$ and $q$ is a prime.
 In the case when $R\leq D$ it is evident.
Now assume that  $R\leq \Phi (G)\cap H$. Then $R\nleq D$ and hence $|R|=p$ for some prime $p$. 
Let $G_{p}$ be a Sylow $p$-subgroup of $H$. Then $G_{p}$ is a Sylow 
$p$-subgroup of $G$  since $H$ is a Hall subgroup of $G$ by Claim (3).
Suppose that  $R\nleq \Phi (H)$. Then for some maximal subgroup $M$ of $H$ 
we have $H=R\rtimes M$, so $G_{p}=R\rtimes (M \cap H)$. But then 
 $R$ has a complement in  $G$ by Gasch\"utz's Theorem 
\cite[I,  17.4]{hupp}. This contradiction shows that $R\leq \Phi (H)$. Suppose that  
  $H/R$
 is cyclic. Then $H$ is nilpotent and so  $\Phi (H)$ is a maximal subgroup of $H$.
 It follows that $H$    is  a cyclic group of order $p^{n}$, where $n > 1$,  a contradiction. 
 Therefore the hypothesis holds for $G/R$. 

If $R\leq \Phi (G)$, then from Lemma 2.8  we deduce that
 $D$ is nilpotent since $DR/R\simeq
D$ is nilpotent, a contradiction. Hence    $\Phi (G)=1$. Therefore $F(D)$ 
is the direct product of some minimal normal subgroups of $G$ by 
Lemma 2.3 since $\Phi (F(G))\leq \Phi (G)=1$.  
 If $G$ has a minimal normal subgroup $L\ne R$ such
 that $L\leq D$, then  $D\simeq D/1=D/ R\cap L$ is nilpotent. Therefore
 $R$ is the unique
 minimal normal subgroup of $G$ contained in $D$. Therefore 
$R=F(D)=C_{D}(R)$. Finally, since $D$ is supersoluble,  a Sylow $p$-subgroup 
$P$ of $G$, where $p$ is the largest prime dividing $|D|$, is normal in 
$D$ and  so $P\leq C_{D}(R)=R$.  Hence  $R$ is the Sylow $p$-subgroup of $D$, so $R$
 is the Sylow $p$-subgroup of $G$ since $D$ is a Hall subgroup of $G$ by Claim (3).

(***) {\sl  $R=C_{G}(R)$. Hence $F(G)=R$.}

Let $C=C_{G}(R)$ and $S$ be a $p$-complement of $C$. Then, in view of Claim (**), 
 $C=R\times S$ is normal in $G$ and so $S$ and $S\cap D$ are  normal in $G$.  
 Therefore Claim (**) implies that $S\cap D=1$. Therefore $S\leq O_{\pi'}(G)=1$ by Claim (6). 
 Hence $C=C_{G}(R)$.

 {\sl Final contradiction for (9).} First assume that $H$ is a $q$-group for some prime $q$
 and 
 $V$ and $W$ are different 
maximal subgroups of $H$. Then $DV$ and $DW$ are supersoluble by Claim (2) and 
$G= DVW=(DV)(DW)$. Hence   $G$ is metanilpotent and then $G/R$ is 
nilpotent by Claim (***). Hence $D=R$ is nilpotent. This contradiction shows 
that $H=AB$, where $A$ is a Sylow $q$-subgroup of $H$ for some prime $q$ dividing $|H|$ and
 $B\ne 1$ is a  $q$-complement of $H$. Let $S$ be a $p$-complement of $D$ such that $SB=BS$.     
Then $DB/F(DB)=DB/R\simeq SB$ is abelian. Hence $|G:C_{G}(S)|$ is a $\{p, q\}$-number. Similarly,
 one can  obtain that $|G:C_{G}(S)|$ is $(\{p\} \cup \{q'\})$-number.
  Hence $|G:C_{G}(S)|$ is a power of $p$.  Therefore a $p$-complement of 
$G$ is supersoluble, which implies that $D=R$, a 
contradiction. Hence we have (9).

(10) {\sl Every chief factor $K/L$ of $G$ below $D$ is non-cyclic.}

If $L\ne 1$, it is true by Claim (1) and the choice of $G$. On the
 other hand, in the case when 
 $L=1$   
  $K$ is  not cyclic by Claim (8).

From Claims (1)--(10) it follows that   Assertions (i)--(vi) are true for $G$, 
which contradicts the choice of $G$. This completes the proof of the 
necessity.

{\sl Sufficiency.}    Let $A$ be a non-identity subgroup of $G$.
 We shall show that $A$ is either 
$\frak{U} $-subnormal or $\frak{U} $-abnormal in $G$. It is clear that $A=V\rtimes W$, where 
$V=A\cap D$ and $W$ is a Hall $\pi'$-subgroup of $H$. Moreover, since $G$ is soluble  and 
 $H$ is a Hall $\pi'$-subgroup of $G$, we can assume without loss
 of generality that  $W\leq H$ and so $A=(A\cap D)(A\cap H)$. If 
$H\leq A$, then $A$ is  $\frak{U} $-abnormal in $G$ 
since $H$ is a Gasch\"utz subgroup     of $G$ by hypothesis. Assume that $A\cap H < H$ and
 let $E=D(A\cap H)$.  Then $E$ is  $\frak{U} $-subnormal in $G$ and $E$ is 
supersoluble by Assertion (vi). Hence $A$ is $\frak{U} $-subnormal in 
$G$ by Lemma 2.1(3).

{\bf Proof of Theorem B.}  {\sl Necessity.}  Suppose that this     is false and
 let $G$ be a counterexample 
 of   minimal order.    
Then $ \Phi\ne D$ and, in view of Theorem A,  
$\Phi \ne 1$.

 (1) {\sl  Assertions (i)--(vi) in Theorem A are true for $G/\Phi$. Moreover,  
Assertions (iii)--(vi) in Theorem A are true for $G/R$ for
 any non-identity normal subgroup $R$ of $G$ not containing $D$. }

 First note that $\Phi (G/\Phi)\cap (G/\Phi)^{\frak{U}}=(\Phi 
(G)/\Phi)\cap (D/\Phi)= (\Phi(G)\cap D)/\Phi= 1$, so  $G/\Phi$ is an 
$E_{\frak{F}}$-group by Lemma 2.1(2). Therefore 
Assertions (i)--(vi) in Theorem A are true for $G/\Phi$.
In order to prove the second assertion of (1), it is enough to
 show that the  hypothesis holds for $G/R$. First note that $G/R\not \in  
\frak{U}$ since $D\nleq R$ by our hypothesis on $R$. 
Let $A/R$ be a subgroup of $G/R$ contaning  
 $\Phi(G/R)\cap (G/R)^{\frak{U}}=\Phi(G/R)\cap (DR/R)$.  Then, since
 $\Phi R/R \leq  \Phi (G/R)$,
 $A$ contains $\Phi$  and so  $A$ is either 
$\frak{U} $-subnormal or $\frak{U} $-abnormal in $G$. Hence $A/R$  is either 
$\frak{U} $-subnormal or $\frak{U} $-abnormal in $G/R$ by Lemma 2.1(2). Therefore
 the  hypothesis holds for $G/R$ for any   normal subgroup $R$ of $G$  not containing $D$. 

(2) {\sl   $G=D\rtimes H$ is a soluble group, where  $H$ is a  
  Hall subgroup of $G$ such that  $H\Phi /\Phi$ is a Gasch\"utz subgroup
 of $G/\Phi$}.

From Claim (1) we get that   
 $D/\Phi =(G/\Phi)^{\frak{U}}$ is a Hall subgroup of $G$.
Therefore $G$ is soluble. 
 It follows also that $\Phi$  is the Sylow $p$-subgroup of $D$ for some prime $p$ and a
 $p$-complement $S$ of $D$ is a Hall subgroup of $G$.
 
Since  $ \Phi  \leq  \Phi (G)$,  from  Lemma 2.8 it follows that $D=R\times S$.
Hence for some minimal normal subgroup $L$ of $G$  we have $L\leq S$ since 
$\Phi \ne D$. Then  $D/L 
=(G/L)^{\frak{U}}$ is a Hall subgroup of $G$ by Claim (1) and hence $\Phi$ is a Sylow 
subgroup of $G$. But this is impossible since $\Phi\leq  
\Phi (G)$.  Hence  $D$ is a Hall subgroup of $G$, so $D$ has a 
complement $H$ in $G$.  It is clear that $H\Phi /\Phi$ is a complement to 
$D/\Phi =(G/\Phi )^{\frak{U}}$ in $G/\Phi$, so $H\Phi /\Phi$ is a Gasch\"utz subgroup
 of $G/\Phi$ since Assertion (i) in Theorem A is true for $G/\Phi $ by 
Claim (1).

(3) {\sl Assertion (iii) in Theorem A is  true for $G$.}

Indeed, from Claim (1) we get that 
  $|(G/\Phi):(D/\Phi)(G'\Phi/\Phi)|=|G:DG'\Phi|=|G:DG'|$ since $\Phi\leq D$.

 (4) {\sl Assertion (iv) in Theorem A is  true for $G$.}

 By Claim (1),   $D/\Phi$ is nilpotent. 
 But then $D$ is nilpotent by Lemma 2.8.

(5) {\sl Assertion (v) in Theorem A is  true for $G$} (This directly follows from Claim (1)).

(6) {\sl Assertion (vi) in Theorem A is  true for $G$.}

  Let $E$ be any  proper subgroup of $G$ containing $D$.  We 
shall show that $E$ is supersoluble.  Suppose that this is false.
Let  $R$ be a minimal 
normal subgroup of $G$ contained in $D$.   
  Then   
 $D/R =(G/R)^{\frak{U}}\leq E/R < G/R$. Hence $E/R$
 is supersoluble by Claim (1). Moreover,   $E/\Phi
 $ is  supersoluble by the same Claim, so  in the case when $H$ is abelian $E$
 is supersoluble by Lemma 2.8. Hence  $H$ is not 
abelian. In this case $D$ is nilpotent by Claim (1).  If  $G$ 
has a  mimimal normal subgroup $L\ne R$ 
such that $L\leq D$, then $E\simeq E/R\cap L$ 
is supersoluble.  Therefore   $R$ is the only minimal 
normal subgroup of $G$ contained in $D$.  Hence $D$ is a Sylow $p$-group of $G$ for 
some prime $p$.   
  Suppose that $R\leq \Phi (D)$. 
Then  $E/\Phi (D)$ is supersoluble.   
 Moreover, $\Phi (D)\leq 
\Phi (E)$ since $D$ is normal in $E$ and hence $E$ is supersoluble. 
Therefore   $\Phi (D)=1$, so  $R=D$ by  the Maschke's Theorem. But then 
$\Phi=1$, a contradiction.  Hence we have (6). 

From Claims (2)--(6) it follows that  the necessity conditions of the theorem   are true for $G$, 
which contradicts the choice of $G$. This completes the proof of the 
necessity.

The sufficiency condition in the theorem directly follows form Theorem A.

\section{Final remarks}

1.  The structure of  ${\mathfrak U}$-critical groups are well-known  
\cite{hu, D}. In particupar, the supersoluble residual of $G^{{\mathfrak 
U}}$  of an ${\mathfrak U}$-critical group  $G$ is a Sylow subgroup of $G$.
  This observation and   Theorem A are    
 motivations for the following question: {\sl Let $G$ be an 
$E_{\frak{U}}$-group. Is  it true then that
 $G^{{\mathfrak U}}$ is a Sylow subgroup of $G$ or, at least, the number 
 $|\pi (G^{{\mathfrak U}})|$   is limited to the top?  }

  The following elementary 
example shows that the answer to this question is negative. 

{\bf Example 4.1. }  
 Let $p_{1} < p_{2} < \cdots < p_{n}  < p$ be a set of primes, $B$ a group of order $p$ 
and $P_{i}$ a simple  ${\mathbb F}_{p_{i}}B$-module
 which is     faithful  for $B$. Let $A_{i}=P_{i}\rtimes B$
 and $G=(\ldots ((A_{1}\Yup A_{2})\Yup A_{3})  \Yup  \cdots )\Yup A_{n}$  (see \cite[p. 50]{hupp}). 
Then $G$ is an  $E_{\frak{U}}$-group,  $G^{\frak{U}}= P_{1}P_{2} \cdots  P_{n} $ 
and $|G/G^{\frak{U}}|=p$.

2.  The following example shows that the subgroup $D$ in Theorem $A$ is not necessary 
nilpotent.

{\bf Example 4.2. } Let $H=H_{2}\rtimes H_{3}$ is a $2$-closed
 Schmidt group, where $H_{2}$ is a Sylow 2-subgroup of $G$ and $H_{3}=\langle a \rangle $
 a cyclic sylow $3$-subgroup  of $G$. Then, by \cite[B, 10.7]{DH}, there exists 
   a simple  ${\mathbb F}_{7}H$-module $P$ 
 which is     faithful  for $H$.   
  Let $G=P\rtimes H$.    It is no difficult to show that  $G$ is an $E_{\frak{U}}$-group and 
$G^{\frak{U}}=PH_{2}$ is non-nilpotent.

3. It is  also  not difficult  to show that  the subgroup $H$ in Theorem A is not 
necessary cyclic.

4.  We do not know the answer to the following question: {\sl What is
 the structure of the group $G$ provided that each nontrivial nilpotent subgroup of $G$
 is either $\frak{U}$-subnormal or  
or $\frak{U}$-abnormal  in $G$?}

5. Partially,  the results of this paper were announced in  \cite{OO}.


\begin{thebibliography}{s2}


\bibitem{Shem} L.A. Shemetkov, \emph{Formations of finite groups},
Moscow, Nauka, Main Editorial Board for Physical and Mathematical
Literature, 1978.



\bibitem{DH}  K. Doerk, T. Hawkes,  \emph{Finite Soluble Groups},
 Walter de Gruyter, Berlin, New York, 1992.




\bibitem{2}  G. Ebert, S. Bauman, A note on subnormal and abnormal chains,
 \emph{J. Algebra}, {\bf 36}(2) (1975), 287--293.




\bibitem{1}  A. Fattahi, Groups  with only normal and abnormal subgroups,
 \emph{J. Algebra}, {\bf 28}(1) (1974), 15--19.


\bibitem{3}   P. F\"orster,  Finite groups all of whose subgroups are  $\frak{F}$-subnormal or  
or $\frak{F}$-subabnormal, \emph{J. Algebra}, {\bf 103}(1) (1986), 285--293.

\bibitem{I}  V.N. Semenchuk,   The structure of finite groups with the $\frak{F}$-abnormal
 or $\frak{F}$-subnormal subgroups, in "Questions of  Algebra",  
 Minsk: Publishing House "University" , {\bf 2}  (1986),  50--55.

\bibitem{II}  V.N. Semenchuk, Finite groups with $f$-abnormal and  или $f$-subnormal subgroups,
 \emph{Mat. Zametki}, {\bf 55}(6) (1994), 111-115. 


\bibitem{III}   V.N. Semenchuk, S.H. Shevchuk,  Finite groups whose  primary subgroups are 
    $\frak{F}$--abnormal    or $\frak{F}$--subnormal, \emph{Izv. Vyssh. Uchebn.
 Zaved. Mathematics,}   {\bf 11} (2011), 46--55. 

\bibitem{Car}  R. Carter, T. Hawkes, The $\frak{F}$-normalizers  of a 
finite soluble group, \emph{J. Algebra}, {\bf 5}(2) (1967), 175--202.


\bibitem{hu} B. Huppert, Normalteiler and maximale Untergruppen
endlicher Gruppen, \emph{Math. Z.}, {\bf 60} (1954), 409--434.


\bibitem{D} K. Doerk, Minimal nicht uberauflosbare, endliche Gruppen,
\emph{Math. Z.}, {\bf 91} (1966),  198--205.


\bibitem{Shem-Sk}  L.A. Shemetkov,  A.N. Skiba,  Formations of Algebraic
Systems,  Nauka, Main Editorial Board for Physical and Mathematical
Literature, Moscow, 1989.




\bibitem{Guo} W. Guo,  \emph{The Theory of Classes of  Groups}, Science
 Press-Kluwer Academic Publishers, Beijing--New York--Dordrecht--Boston--London, 2000.





\bibitem{BE} A. Ballester-Bolinches, L.M. Ezquerro, \emph{Classes of
Finite groups},  Springer, Dordrecht, 2006.




  \bibitem{5} W. Guo, A.N. Skiba,
 On $\frak{F}\phi ^*$-hypercentral subgroups of
finite groups, \emph{J. Algebra}, {\bf 372} (2012), 275--292.


\bibitem{hupp} B. Huppert, \emph{Endliche Gruppen I},
Springer-Verlag, Berlin, Heidelberg, New York, 1967.


\bibitem{We} M. Weinstein (ed.), et al.,  \emph{Between Nilpotent and Solvable}, Polygonal
Publishing House, Passaic N. J., 1982.


\bibitem{11}  A.F.Vasil'ev, T.I. Vasil'eva,
V.N.Tyutyanov, On the finite groups of supersoluble type, \emph{Sib. Math. J.},
{\bf 51}(6) (2010),   1004--1012.

\bibitem{22} A.F.Vasil'ev,
T.I. Vasil'eva, V.N.Tyutyanov, On the rpoducts of $\mathbb P$-subnormal subgroups of finite groups,
\emph{Sib. Math. J.}, {\bf 51}(6) (2012),   47--54.


\bibitem{HupBl} B. Huppert, N. Blackburn, \emph{Finite Groups II},
Springer-Verlag, Berlin, New-York, 1982.


\bibitem{OO}   V.N. Semenchuk,    Finite groups with generalized subnormal
 formation subgroups,  \emph{Proc. Fr.  Skorina Gomel 
St.  Univ.},  {\bf 84}(3) (2014),  104--107.




\end{thebibliography}
\end{document}